\documentstyle[12pt]{article}
\begin{document}
\title{Singular Measures in  Circle Dynamics}
\author{Jacek Graczyk  \\
Institute of Mathematics, Warsaw University, \\
ul.Banacha 2, 02-097 Warszawa, Poland.
\and
Grzegorz \'{S}wi\c{a}tek
\thanks{Partially supported by NSF grant \#431-3604A} \\
Institute for Mathematical Sciences,\\
SUNY at Stony Brook,\\
Stony Brook, NY 11794, USA.}
\date{July 23, 1992}
\maketitle
\newtheorem{lem}{Lemma}[section]
\newtheorem{conjec}{Conjecture}
\newtheorem{coro}{Corollary}[section]
\newtheorem{prop}{Proposition}
\newtheorem{con}{Construction}[section]
\newtheorem{defi}{Definition}[section]
\newcommand{\hf}{\hat{f}}
\newcommand{\B}{{\cal B}}
\newcommand{\de}{{\bf \delta}}
\newtheorem{fact}{Fact}[section]
\newenvironment{proof}
{{\bf Proof:}\newline}{\begin{flushright}$\Box$\end{flushright}}
\newcommand{\ex}{{\cal EX}}
\newcommand{\Cr}{{\bf Cr}}
\newcommand{\dist}{\mbox{dist}}
\newcommand{\Bo}{\Box^{n}_{i}}
\newcommand{\E}{{\bf E}} 
\newcommand{\Po}{{\bf Poin}}
\newcommand{\DPo}{{\bf DPoin}}
\newcommand{\td}{{\underline \tau}}
\newcommand{\tg}{{\overline \tau}}
\newcommand{\gd}{{\underline \gamma}}
\newcommand{\gu}{{\overline \gamma}}
\newenvironment{double}{\renewcommand{\baselinestretch}{2}\protect\Large
\protect\normalsize}{}

\begin{abstract}
Critical circle homeomorphisms have an invariant measure 
totally singular with respect to the Lebesgue measure. We prove
that singularities  of the invariant measure are of H\H{o}lder type.
The Hausdorff dimension of the invariant measure is less than $1$
but greater than $0$. 
\end{abstract}

\section{Preliminaries}
\subsection{Discussion of the Results}
The long time behavior of nonlinear dynamical systems can be often
characterized by means of invariant measures. A variety of ``multifractal
formalisms'' have been developed recently to study statistical properties
of singular measures (see \cite{HJKPS}, \cite{CJ} for more details) 
which appear as a natural description of many physical phenomena.
One of the characteristic quantities describing the multifractal structure
of a singular measure $\mu$ is a singularity spectrum $g(\alpha)$ which
is usually  defined in an informal way (see \cite{HJKPS},
\cite{CJ} and many others) as follows:
\newline
Cover the support  of $\mu$ by small boxes $L_{i}$ of size $l$. Then
define the 
singularity strength $\alpha_{i}$ of $\mu$ in the $i$-th box by the relation:
\[\mu(L_{L}) \sim l^{\alpha_{i}}.\]
We count the number of boxes $N(\alpha)$ where $\mu$ has singularity strength
between $\alpha$ and $\alpha + d\alpha$ (whatever that is to mean).
Then $g(\alpha)$ is defined by the requirement that
\[N(\alpha) \sim l^{g(\alpha)}.\]

Unfortunately, many ``multifractal formalisms'' suffer from mathematical
ambiguities (see \cite{CJ} for a fuller discussion of this problem;
for example, is $g(\alpha)$ a Hausdorff or a box dimension or something else?)
 even if they provide qualitative information on a given dynamical system.
In the present paper we would like to propose a method of describing
the dynamics of critical circle homeomorphisms. Our method is
more general then the method  relying on the scalings exponents (see \cite{F}),
and on the other hand, mathematically rigorous unlike the 
``multifractal formalism'' in its present shape. 

\paragraph{Description of the method.}
Unlike typical smooth diffeomorphisms, which were treated
in~\cite{HH}, all critical circle
homeomorphisms have singular invariant measures.
Moreover, it turns out that the unique normalized invariant measure is always
completely  singular with respect to the Lebesgue measure.
We introduce two singularity exponents, the lower and the upper one,
to measure the increments of distribution of the invariant measure in 
the logarithmic scale. We study these exponents with respect to two natural
measures on the circle: the invariant measure $\mu$ and the Lebesgue measure
$\lambda$. By ergodicity, these exponents are constants on sets of
full measure $\mu$ or $\lambda$, respectively.

Our main achievement is to prove uniform bounds for the exponents
in the class of circle maps with a critical point of polynomial type
 and an irrational rotation number of constant type\footnote{Constant
type irrational number means
that coefficients in the continued fraction representation are bounded.}

\paragraph{Universality.}
We should mention here that for critical maps with all critical points
of polynomial type and
rotation numbers of algebraical degree $2$, the {\em universality
conjecture} implies that the upper and the lower exponents coincide.
The reader may consult~\cite{R} for more information about circle map 
universality and its consequences. 
There are strong computer-based arguments in favor of the conjecture 
(see \cite{KO}, also for the list of other references). However, in
the absence of a definite rigorous proof,
we continue to regard the conjecture as just that, and will refrain
from using it in our discussion.

Another important quantity which describes a singular measure $\mu$ is 
the Hausdorff dimension  $HD(\mu)$ of the measure theoretical support
(i.e., the infimum of the dimensions of the sets of the full measure).
Using the singularity exponents we immediately obtain universal bounds
on $HD(\mu)$ in our class of circle maps.

\subparagraph{Hausdorff dimension.} The renormalization
group analysis applied to
study high iterates of circle maps with special rotation numbers (like
the golden mean) lead to several universality conjectures (see for example 
\cite{HJKPS}, \cite{KO}, \cite{R}).
  
We state one which is certainly
true provided the {\em golden mean universality conjecture} holds.
\begin{conjec}
$HD(\mu)$ is constant
in any topological conjugacy class of cubic critical homeomorphisms with 
rotation number of algebraical degree $2$.
\end{conjec}
  
An intriguing question remains about universal properties for
more general irrationals. We think that the same conjecture should be
true for any irrational rotation number, even of Liouville type.
However, the evidence for that is scarce and we leave this merely as
an interesting open question.

\subsection{Introduction}
\paragraph{Assumptions.}
All results in this paper are true for $C^{3}$ smooth circle homeomorphisms
with finitely many critical points of polynomial type and 
an irrational rotation number of constant type.

For simplicity of our presentation we will give detailed proofs of our 
results only for maps with exactly one critical point which after a $C^{2}$
change of coordinate system can be written in the proximity of a critical
point $0$ in the form $x_{i} \longmapsto (x-x_{i})^{3} + \varepsilon$.
As a consequence, the circle can be covered by two overlapping sets:
in the vicinity of the critical point $0$  by an symmetric interval
$U$ and a ``remote" interval $V$ on which the first   
derivative is bounded away from zero. On the interval
$V$ the map has strictly negative Schwarzian derivative.
We reserve the letter $f$ for maps from the class defined above. 
The real line is projected to the unit circle by means of the map
\[ x \longmapsto \exp(2\pi x i).\]
Denote by $|x-y|$ the distance between points $x$ and $y$ on the circle
in the metric induced by the projection.
    
\paragraph{Uniform Constants.}
Following the convention of \cite{SW1} we will mean by a uniform constant
a function on our class of maps which continuously depends on
the quasisymmetric norm of the map,
the logarithm the size of $U$, the lower bound of the derivative on the
remote arc and the $C^{3}$ norm. Uniform constants will be always denoted
by the letter $K$. Whenever confusion can arise 
we specify uniform constants by adding subscripts.

\paragraph{Continued Fractions and Dynamics.}
Let $p_{n}/q_{n}$ be the $ n$-th continued fraction approximant of the
rotation number $\rho $ of $f$.
 The numbers $q_{n}$ and the coefficients $a_{n}$ 
in the continued fraction representation of $\rho$ are related  by
 the recurrence formula:
\[ q_{n+1} = a_{n}q_{n} + q_{n-1}, \hspace{5 mm} n\geq 2,\:\:q_{0}=1,\,\,
q_{1}=a_{1}\]
Dynamically $q_{n}$ is that iterate of the rotation by $\rho$ for which
the orbit of any point makes the closest return so far to the point itself.
According to the Yoccoz Theorem (see \cite{Y}) a homeomorphism from our class
 is conjugated to a rotation. In particular, it implies the same order
 of orbits both for $f$ and the rotation by $\rho$.
The numbers $q_{n}$ are called {\em closest returns}.

\paragraph{Continued Fractions and Partitions.}
We will use the orbit of a critical point $0$ to define a system of partitions
of the circle. First, we define two sets of closed intervals of order $n$: 
\begin{tabbing}
\= \hspace{10 mm} \= $q_{n-1}$  ``short" intervals:
 \= \hspace{ 10. mm} \= $(z,f^{q_{n}}(z))$,

\ldots,$f^{q_{n-1}-1}(z,f^{q_{n}}(z))$.\\
\>and \> \> \>  \\
\> \>$q_{n}$  ``lengthy" intervals: \>  \>$ (z,f^{q_{n-1}}(z))$,
\dots,$f^{q_{n}-1}(z,f^{q_{n-1}}(z))$.\\
\end{tabbing}
The  ``lengthy" and ``short" intervals are mutually disjoint except for 
the endpoints and cover the whole circle.
The partition obtained by the above construction 
will be denoted by ${\cal B}(n;f)$ and called the dynamical partition of
the $n$-th order.

We will briefly explain the structure of the dynamical partitions.
Take two subsequent dynamical partition of order $n$ and $n+1$.
The latter is clearly a refinement of the former.
All ``short" intervals of ${\cal B}(n;f)$ become the ``lengthy" intervals 
of $\B(n+1;f)$ while all ``lengthy" intervals of $\B(n)$ are split into
$a_{n}$ ``lengthy" intervals and $1$ ``short" interval of the next partition
$\B(n+1;f)$.
An interval of the $n$-th dynamical partition
will be denoted by $\Box^{n}(f)$ or by  $\Box^{n}_{x}(f)$ if we want
to emphasize that the interval contains a given point $x$.

We will drop $f$ in the the notation when no confusion can arise. 

\paragraph{Bounded Geometry.}
Let us quote a few basic
results about the geometry of dynamical partitions which are commonly
 referred to as ``bounded geometry".(see for the proof \cite{H} and \cite{SW1}) 
\begin{itemize}
 \item The ratio of lengths of two adjacent elements of any dynamical partition
is bounded by a uniform constant.
\item For any element of any dynamical partition, the ratios of its length to
the lengths of extreme intervals of the next partition subdividing it are 
bounded by a uniform constant.
\end{itemize}

As a corollary we obtain that the elements of the $n$-th dynamical
partition are exponentially small.
\begin{fact}\label{1}
There are uniform constants $K_{1}$, $K_{2} \leq 1 \,$,$K_{3} \leq 1$ so that
     \[ K_{1} K_{2}^{n}\leq |\Box^{n}|\leq K_{1}K_{3}^{n}\]
holds for all natural numbers $n$.
\end{fact}

\section{Technical Tools}

\paragraph{Distortion Lemma.}

We will call a chain of intervals a sequence of intervals such that each is
mapped onto the next by the map $f$.

Denote by $\Cr(a,b,c,d)$ a cross-ratio  of the quadruple $(a,b,c,d)$,
$a < b < c < d$  given by the formula
\[ \Cr(a,b,c,d) = \frac{|b-a||d-c|}{|c-a||d-b|}. \]
Here is one of possible stating the Distortion Lemma for critical circle 
homeomorphisms:
\begin{lem}\label{DL}
Take  a chain of disjoint intervals 
\[ (a_{0},b_{0}),\dots, (a_{m},b_{m})\]
which do not contain a critical point of $f$.
Then, for arbitrary points $x,y \in (a_{1},b_{1})$, the uniform estimate
\[ \big|\log\frac{(f^{m})'(x)}{(f^{m})'(y)}\big|
 \leq K \Cr(f^{m}(a_{0}),f^{m}(x),f^{m}(y),f^{m}(b_{0}))\] 
holds.
\end{lem}

\paragraph{The Pure Singularity Property.}
To have a ``dynamical measure" of size of an interval
we will make the following definition:
\begin{defi}
An interval $J$ will be said of the $j$-th order of size if
\[ j  =\max\{i:\forall_{x\in J} f^{q_{i}}\not\in J\} +1. \]
\end{defi}
Note that each interval of a chain is of the same order of size.

Let us introduce  a one form
\[{\cal N}f = \frac{f''}{f'}dx\]
called a nonlinearity of $f$. As opposed to diffeomorphisms,
the nonlinearity of critical circle maps which measures the distortion
on chains of disjoint intervals, is non-integrable.

One of the main achievements of \cite{SW2} was  that the distortion
coming from parts of the circle far away from critical points can be
neglected with an almost exponentially small error with the order of
size of a given chain. It means that asymptotically only what happens
in the small neighborhood of a critical point matters. 

We pass to a detailed formulation of the {\bf Pure Singularity Property}. 
Suppose we have a chain of disjoint intervals
\[ (a_{0},b_{0}),\ldots, (a_{m},b_{m})\]
of the k-th order of size and symmetric neighborhood $U_{j}$ with size
of the order $j$.
Then
\[ |\int_{C_{j}} {\cal N}f| \leq K\exp(-\sqrt{k-j}),\]
where $C_{j}$ is a union of these intervals of the chain which are not
contained in $U_{j}$ and a constant $K$ is uniform.

\paragraph{Integral Formula.}
We introduce another cross-ratio $\Po (a,b,c,d)$ of a given quadruple
$(a,b,c,d)\,$, $a \leq b \leq c \leq d$, by the following formula:
\[ \Po(a,b,c,d) = \frac{|b-c||d-a|}{|c-a||d-b|}. \]
 By the distortion of the cross-ratio $ \Po(a,b,c,d)$ by $f$ we mean
\[\DPo(a,b,c,d;f) = \frac{\Po(f(a),f(b),f(c),f(d))}{\Po(a,b,c,d)}. \]
There is a very simple relation between cross-ratios $\Cr$ and $\Po$. Namely,
\[\Cr = \frac{1}{1+\Po}. \]
The logarithm of the distortion of the cross-ratio $\Po$ can be expressed
by the integral formula. The formula is due to Sullivan \cite{S}:
\[ - \log(\Po(a,b,c,d)) = \int\int_{S}\frac{dxdy}{(x-y)^{2}},\]  
where $S = \{(x,y)\,\, a \leq x \leq b \mbox{ and } c \leq y \leq d \}$.

Consequently,
\[\DPo(a,b,c,d;f) = \int\int_{S}d\mu-(f^{*}\times f^{*})d\mu,\]
where $\mu$ is equal to $\frac{dxdy}{(x-y)^{2}}$.
Calculating the integrand we get that
\[d\mu-(f^{*}\times f^{*})d\mu = \Big( 1 - \frac{f'(x)f'(y)}{(\frac{f(x)-f(y)}
{(x-y)})^{2}}\Big)\frac{dxdy}{(x-y)^{2}}.\]

For maps with negative Schwarzian derivative the integrand is positive and,
as a consequence, the cross-ratio is not decreased by $f$.
In the next paragraph we estimate how much the cross-ratio $\Po$ is
expanded by maps with strictly negative Schwarzian.

\paragraph{Expansion Lemma.}
Let $a<b<c<d$.
Suppose we have a chain of disjoint intervals
\[ (a_{0},b_{0}),\dots, (a_{m},b_{m})\]
of the n-th order of size which omit a critical point $0$.
Then
\[\log(\DPo(a,b,c,d;f^{m})) \geq \] 
\[  K_{1}\sum_{f^{i}(a,d)\subset U}
\frac{|f^{i}(a)-f^{j}(b)||f^{i}(c)-f^{j}(d)|}
{\max(|f^{i}(a)|,|f^{i}(d)|)^{2}} + K_{2}\exp(-\sqrt{n}),\]
where $K_{1}$ and $K_{2}$ are uniform constants.

\begin{proof}
By the Pure Singularity Property we get that
\[\sum_{V}\log(\DPo(a_{i},b_{i},c_{i},d_{i};f))
\leq K_{3}\exp(-\sqrt{n}),\]
where the sum $\sum_{V}$ is over all these indexes $i$ for which 
$f^{i}(a,d)$ intersects $V$.
Next we use Integral Formula to estimate the expansion of the cross-ratio
for quadruples $(a_{i},b_{i},c_{i},d_{i})$ contained in $U$.

\begin{equation}\label{equ:15na,1}
\log(\DPo(a_{i},b_{i},c_{i},d_{i};f)) =
\int\int_{S_{i}}
\Big(1-\frac{9x^{2}y^{2}}{(\frac{x^{3}-y^{3}}{x-y})^{2}}\Big)
\frac{dxdy}{(x-y)^{2}} \; .
\end{equation}

Here, $S_{i}$ is defined by:
$S_{i}=\{(x,y):\,\, f^{i}(a) \leq x \leq f^{i}(b),\,\,\,
 f^{i}(c) \leq y \leq f^{i}(d)\}$ .

By algebra, the
right-hand of equation~\ref{equ:15na,1} is rewritten as:

\begin{equation}\label{equ:15na,2}
\int\int_{S_{i}}\frac{(x^{2}+xy+y^{2})^{2}-(3xy)^{2}}{(x^{2}+xy+y^{2})^{2}}
\frac{dxdy}{(x-y)^{2}} =
\end{equation} 
\[=\int\int_{S_{i}}\frac{x^{2}+4xy+y^{2}}{(x^{2}+xy+y^{2})^{2}}dxdy \geq
 \frac{|f^{i}(a)-f^{j}(b)||f^{i}(c)-f^{j}(d)|}
{3\max(|f^{i}(a)|,|f^{i}(d)|)^{2}} \]
which immediately gives the claim of the Expansion Lemma. The last
inequality follows if we forget the numerator while dropping the
power of the denominator by $1$, and next estimate the
denominator by  
\[ 3 \max (|f^{i}(a)|,|f^{i}(d)|)^{2}\; .\]
\end{proof}

\section{ Singularity of the Invariant Measure}
It is a well known fact that homeomorphisms of the circle have exactly
one invariant measure $\mu$.  In this section
we will investigate the properties of this measure for critical circle
homeomorphisms.
We will start with the following observation.
\begin{prop}
The invariant measure $\mu$ is totally singular with respect
to the Lebesgue measure.
\end{prop}
\begin{proof}
Let $\phi$ be the conjugacy between $f$ and a rotation $\rho$,
$ \rho \circ \phi = \phi \circ f $. It is enough to show that $\phi$
 has the first derivative equal to zero on a set of full Lebesgue measure.
To the contrary, suppose that at some point $x$ the first derivative exists
and is non-zero. Consider a first return $q_{n}$. The $q_{n+1} - 1$ images
of $(x,f^{-q_{n}}(x))$ are disjoint. Clearly, there is an infinite
sequence of first returns so that $f^{q_{n+1}}$ on this interval is not
a diffeomorphism. By our conjugacy assumption, this map must be arbitrary
$C^{0}$ close to a linear map for large values of $n$. On the other hand, 
by bounded geometry, it is a composition of a few bounded distortion
diffeomorphisms and a bounded number of critical iterates which are not 
diffeomorphisms. But maps of this type can not be arbitrarily $C^{0}$ close
to linear.
\end{proof}

Another important property is ergodicity.

\begin{prop}\label{ergodic}
The map $f$ is ergodic with respect to the Lebesgue measure
 $\lambda$.
\end{prop}

\begin{proof}
Suppose that there exist an invariant set $A$ of positive but not full 
the Lebesgue measure $\lambda(A)$.
 
We fix $\varepsilon > 0$. Then by the Lebesgue Density Theorem we
can find a point $z$ and a number $n_{0}$ so that for all $n \geq n_{0}$
the Lebesgue measure of an interval of $n$-th partition which contains $z$
satisfies the inequality
\[\lambda(\Box^{n}_{z}\cap A) \geq (1-\varepsilon) |\Box^{n}_{z}|\]
or, equivalently, 
\[\lambda(\Box^{n}_{z}\cap A^{c}) \leq \varepsilon  |\Box^{n}_{z}|,\]
where $A^{c}$ denotes the compliment of $A$.

Taking $q_{n+1}+q_{n}$ or $q_{n}+q_{n-1}$ images of $\Box^{n}_{z}$
in dependence
on  $\Box^{n}_{z}$ is a ``short" or a ``long" interval of the $n$-th dynamical
partition we obtain a cover of the circle. One can check that each point
of the circle belongs to at most two intervals of this cover.
We want to estimate $\lambda(f^{k}( \Box^{n}_{z})\cap f^{k}(A^{c}))$ for each
interval of the cover.

If  $f^{i}(\Box^{n}_{z})$ contains a critical point then there
is a uniform constant $K_{1}$ so that

\[ \frac{\lambda(f^{i+1}( \Box^{n}_{z} \cap A^{c}))}{|f^{i+1}( \Box^{n}_{z})|}
\leq K_{1}
 \frac{\lambda(f^{i}( \Box^{n}_{z} \cap A^{c}))}{|(f^{i}( \Box^{n}_{z})|}\]
The above inequality and the Distortion Lemma implies  that
\[\frac{\lambda(f^{k}( \Box^{n}_{z})\cap f^{k}(A^{c}))}
{|f^{k}( \Box^{n}_{z})|} \leq K_{2}
\frac{\lambda( \Box^{n}_{z}\cap A^{c})}{|( \Box^{n}_{z})|}.\]
Since $A^{c}$ is invariant we obtain that

\[ \lambda(A^{c}) \,\, \leq \,\,
 \sum_{k} \lambda(f^{k}( \Box^{n}_{z})\cap f^{k}(A^{c})) \,\, \leq
 \hspace*{\fill}\]
\[ K_{2} \sum_{k} |f^{k}( \Box^{n}_{z})|
 \frac{\lambda( \Box^{n}_{z} \cap A^{c})}
{|\Box^{n}_{z}|} \,\, \leq \,\, K_{2}\varepsilon,\]
in contradiction to our assumption that $\lambda(A^{c})$ is positive.
\end{proof}

\paragraph{Singularity exponents.}
We are going to study the nature of singularities of an invariant measure 
${ \mu}$ using some ideas underlying  the concept
of multifractal measures and multifractals, the objects which are intensively
studied by physicists.
Let us discuss briefly the concept of
 a singularity exponent of an invariant measure which can be loosely
defined in the following way: 
Let $M(x) = \int_{0}^{x} d \mu$ be the distribution function of 
measure $\mu$. 
If the increments in $M(x)$ between two close points $x$ and $x+ \epsilon$
are of the order $\epsilon^{\tau(x)}$ then we will say that the distribution
$M(x)$ has in the point $x$ an exponent of singularity  $\tau(x)$.

For mathematical exactness we will introduce
two exponents of singularity, the upper and the lower one.

\begin{defi}
Let $\mu$ be a measure completely singular with respect to $\lambda$
with distribution function $M(x)$.
Then by the upper and the lower singularity exponents we mean respectively
\[{\overline \tau(x)} = \limsup_{\epsilon \rightarrow 0}
\frac{\log(M(x+\epsilon)-M(x))}{\log(|\epsilon|)}\]
and
\[{\underline \tau(x)} = \liminf_{\epsilon \rightarrow 0} 
\frac{\log(M(x+\epsilon)-M(x))}{\log(|\epsilon|)}.\]
\end{defi}
Taking into consideration that the Lebesgue measure is the image of
 $\mu$ by the conjugating homeomorphism $\phi$, we can rewrite
 the exponents $\td$ and $\tg$ in the language of
the dynamical partitions \footnote{ Here we use the fact
that the rotation number $\rho$ is of bounded type.}
\[\tg(x) =\limsup_{n \rightarrow \infty} \frac{\log|\phi(\Box^{n}_{x})|}
{\log|\Box^{n}_{x}|}\]
and
\[\td(x) =\liminf_{n \rightarrow \infty} \frac{\log|\phi(\Box^{n}_{x})|}
{\log|\Box^{n}_{x}|}\].
\paragraph{The exponents are constants.} 
The Distortion Lemma immediately implies that 
\begin{lem}
The exponents $\td(x)$ and $\tg(x)$ are f invariant.
\end{lem}
By Proposition~\ref{ergodic} and the uniqueness of the invariant measure
$\mu$ we get that
\begin{itemize}
\item For almost all points with respect to the Lebesgue measure
the exponents are constants. We will denote these constants by
 $\td(\lambda)$ and $\tg(\lambda)$ respectively.
\item The above statement holds verbatim 
if ``the Lebesgue measure" is replaced by $\mu$. Denote these new constants by
 $\td(\mu)$ and $\tg(\mu)$ respectively. 
\end{itemize}

We pass to the formulation of our Main Theorem.
\paragraph{The Main Theorem.}
The singularities of the invariant measure $\mu$ are of H\H{o}lder type.
It means that there exist uniform constants $K_{1}$ and $K_{2}$ so that
for almost all $x$ in the sense of the measure $\mu$ the following estimates
\[ 0 < K_{1} < \td(\mu) \leq \tg(\mu)< K_{2} <1 \]
hold.
\paragraph{Remark.}
We should mention here that $\td(\lambda)$ and $\tg(\lambda)$
are uniformly greater than  $1$ and less than infinity.

The proof of the Main theorem will occupy the whole next section.
 
\paragraph{Reformulation of the Main Theorem.}
For technical reasons we introduce new exponents $\gd(x)$ and $\gu(x)$ 
which lives in the phase space of the rotation $\rho$
\[\gd(x)= \tg^{ -1}(\phi^{-1}(x)) \hspace{3 mm}{\mbox and}
\hspace{3 mm} \gu(x)= \td^{-1}(\phi^{-1}(x))\]
and state the Main Theorem in the following equivalent form:
\paragraph{}
There are uniform constants $K_{1}$ and $K_{2}$ so that 
for almost all points $x$ with respect to the Lebesgue measure $\lambda$ the
estimates 
\[ 1 < K_{1} < \gd(x) \leq \gu(x)< K_{2} <\infty\]
hold.

\section{Proof of the Main Theorem}
\subsection{Discrepancy} 
Our main object in this paragraph is to establish a quantity which would
measure  nonlinear behavior of critical maps. We want to show that critical
maps stay away in a certain uniform distance from diffeomorphisms.
To this end we will introduce a notion of discrepancy.
\paragraph{Discrepancy between partitions.}

We always assume that the length of the interval being partitioned is
less than $1$.
\begin{defi}\label{defi:1,1}
A {\em partition} of $I$, denoted with $P_{I}$, is a set (possibly
infinite) of closed subintervals of $I$, disjoint except for the
endpoints, whose union is $I$. In addition, we assume that the entropy
$H(P_{I})$ is finite.
\end{defi}

Given $J \subset I$. Partition $P_{I}$ induces in natural way a partition
of $J$ denoted by ${\bf[P_{I}:J]}$.

There is a probabilistic measure on $P_{I}$ defined by
\[ \mu(X) := \sum_{w\in X} \frac{|w|}{|I|} \] for every $X\subset
P_{I}$,  where
$|\cdot|$ stands for the Lebesgue measure. 

Two partitions, $P_{J}$ and $P_{I}$ will be considered isomorphic if
there is a homeomorphism $h$ from $I$ to $J$ which maps each element of
$P_{I}$ onto an element of $P_{J}$.

\begin{defi}\label{defi:1,2}
The {\em discrepancy} between isomorphic partitions $P_{I}$ and $P_{J}$
is defined as $\int_{P_{I}}\log_{+}\frac{dh}{d\mu}$  where
$\frac{dh}{d\mu}$ is the Jacobian of the isomorphism, while $\log_{+}$
means $\max(0,\log)$.
\end{defi}
The reader may note that $\de(P_{I},P_{J})$ cannot
be arbitrarily large regardless of the partitions involved.

\paragraph{The Discrepancy Lemma.}
For any $n$ and $r$ the partitions
\[[\B((n+1)r;f):\Box^{nr}(f)]
 \stackrel{\phi}{\longmapsto} [\B((n+1)r;\rho):\Box^{nr}(\rho)]\]  
are isomorphic and the isomorphism is given by the  conjugation $\phi$. 
As it turns out the discrepancy between these partitions
is uniformly bounded away from zero.
\begin{lem}
We can choose $r$ so that the inequality
\[ \de([\B((n+1)r;f):\Box^{nr}(f)], [\B((n+1)r;\rho):\Box^{nr}(\rho)]
 \geq K\]  
is satisfied for large $n$ and a uniform constant $K$.
\end{lem}
\begin{proof}
The interval $\Box^{nr}(f)$ contains at most two critical points of the map
$f^{q_{nr}}$. Bounded Geometry implies that we can choose a number $r$ in
the definition of the refined dynamical partition $[\B((n+1)r;f)]$ so that:
\begin{itemize}
\item There exist three consecutive elements 
\[(a,b)\,,(b,c)\,,(c,d)\,\,\, {\mbox of}\,\,\, [\B((n+1)r;f):~\Box^{nr}(f)] \]
which do not contain a critical point of $f^{q_{nr}}$ and the length of
the interval $(a,d)$ is at least comparable to
$|\Box^{nr}(f)|$, i.e. 
\[ |(a,d)| > K |\Box^{nr}(f)| \] with uniform $K$. 
\item The intervals $f^{q_{nr}}((a,b))\,$, $f^{q_{nr}}((b,c))\,$, 
 $f^{q_{nr}}((c,d))\,$ belong to the partition  $[\B((n+1)r;f):\Box^{nr}(f)] $.
\end{itemize}  
 From the Expansion Lemma we have that 
\[\log(\DPo(a,b,c,d;f^{q_{n}})) \geq K_{1} \frac{|f^{i}(a)-f^{i}(d)|^{2}}
{|f^{i}(d)|^{2}} + K_{2}\exp(-\sqrt{n}),\]
where $f^{i}((a,b))$ is the closest interval to $0$ amongst all $q_{nr}$ images
of $(a,d)$ by $f$. Therefore,
the distortion of the cross-ratio $\Po(a,b,c,d)$ by $f^{q_{nr}}$
is by a definite amount greater than $1$ since $r$ which controls the relative
seize of the elements of $[\B((n+1)r;f):\Box^{nr}(f)]$ is not too large.

But the distortion of the cross-ratio $\Po(\phi(a),\phi(b),\phi(c),\phi(d))$ 
by any iterate of $\rho$ is equal to $1$ since $\rho$ is an isometry.
Hence, the discrepancy between partitions under consideration must be uniformly
separated from zero, provided $n$ is large enough. This concludes the proof.
\end{proof} 
\subsection{Partition Lemma}

Here, we have a lemma about partitions:
\begin{prop}\label{prop:1,1}

Consider intervals $I$ and $J$ with  isomorphic partitions $P_{I}$ and $P_{J}$
respectively. Assume the following:
\[ \int \frac{|\log\mu(h(w))|}{|\log|J||}|\log\mu(w)| d\mu(w) 
 \leq K_{3}\de^{2}(P_{I},P_{J})\; .\]  
If  
\[ \frac{|\log|J||}{|\log|I||} \leq \min(2, 1 +
\frac{K_{1}}{H(P_{I})}  \de^{2}(P_{I},P_{J}))\; , \]
then
\[\sum_{w\in P_{I}} \frac{|\log|h(w)||}{|\log|w||}\mu(w) >
\frac{|\log|J||}{|\log|I||}(1+
K_{2} \frac{\de^{2}(P_{I},P_{J})}{|\log|I||})\; .\]
\end{prop}

We will first work to approximate the sum
\begin{equation}\label{equ:2,1}
\sum_{w\in P_{I}} \frac{|\log|h(w)||}{|\log|w||}\mu(w)
\end{equation}
by a sum easier to deal with. Let us consider an individual term:
\[ \frac{|\log|h(w)||}{|\log|w||}\mu(w) =
\frac{|\log|J||}{|\log|I||}\mu(w)
\frac{1+\frac{|\log\mu(h(w))|}{|\log|J||}}{1+\frac{|\log\mu(w)|}{|\log|I||}}
\; .\]
Now, an expression of the type 
\[\frac{1+x}{1+y}\] 
for positive $x,y$ can be approximated with $1+x-y$ so that 
\begin{equation}\label{equ:2,2}
1+x-y = \frac{1+x}{1+y} + \frac{y(x-y)}{1+y} > \frac{1+x}{1+y} + yx\; .
\end{equation}

Inequality ~\ref{equ:2,2} allows us to bound
a term of sum ~\ref{equ:2,1} from below by
\[ \frac{|\log|J||}{|\log|I||}\mu(w) (1 +
\frac{|\log\mu(h(w))|}{|\log|J||} -\frac{\mu(w)}{|\log|I||} + Q) \]
where the ``quadratic correction'' $Q$ equals
\[ \frac{|\log\mu(h(w))|}{|\log|J||}\frac{|\log\mu(w)|}{|\log|I||}\;
.\]

Let us now bound the contribution of all quadratic corrections to sum
~\ref{equ:2,1}. It is equal to 
\[ \sum_{w\in P_{I}}\frac{|\log|J||}{|\log|I||}\mu(w)
\frac{|\log\mu(h(w))|}{|\log|J||}\frac{|\log\mu(w)|}{|\log|I||}\; .\]

Now we use the first assumption of the proposition to see that this
quantity is than not greater than
\[ \frac{|\log|J||}{|\log|I||}K_{3}\frac{\de^{2}(P_{I}, P_{J})}
{|\log|I||} \; .\]

We can see that to prove Proposition~\ref{prop:1,1} it is sufficient
to show that 
\begin{equation}\label{equ:2,3}
\sum_{w\in P_{I}} (\frac{|\log\mu(h(w))|}{|\log|J||}
-\frac{|\log\mu(w)|}{|\log|I||})\mu(w) > K_{4}\frac{\de^{2}(P_{I}, P_{J})}
{|\log|I||}\; ,  
\end{equation}
that is, to neglect the quadratic corrections. Indeed, we will just
need to pick $K_{3} := K_{4}/2$ to ensure that the quadratic
corrections will not spoil the estimate.

We claim that estimate~\ref{equ:2,3} follows from the following:
\begin{equation}\label{equ:4,1}
\sum_{w\in P_{I}}|\log\mu(h(w))|\mu(w) - \sum_{w\in
P_{I}}|\log\mu(w)|\mu(w) \geq K_{5}\de^{2}(P_{I},P_{J})\; .
\end{equation}

Indeed, assume that~\ref{equ:4,1} holds. The left-hand side of
estimate~\ref{equ:2,3} is
\begin{equation}\label{equ:6,1}
 \sum_{w\in P_{I}} (\frac{|\log\mu(h(w))|}{|\log|J||}
-\frac{|\log\mu(w)|}{|\log|I||})\mu(w) =
\end{equation} 
\[ = \frac{1}{|\log|I||}
 (\frac{|\log|I||}{|\log|J||}\sum_{w\in
P_{I}} |\log\mu(h(w))| - \sum_{w\in P_{I}}|\log\mu(w)|)\; . \]

We know by hypotheses of Proposition~\ref{prop:1,1} that 
\[ \frac{|\log|J||}{|\log|I||} = 1 +
K_{6}\frac{\de^{2}(P_{I},P_{J})}{H(P_{I})} \] where $K_{6}$ is not
greater than a certain constant $K_{1}$ which we are free to specify,
and, in addition, this quantity is not greater than $2$.

From this and estimate~\ref{equ:4,1} we can bound expression~\ref{equ:6,1}
from below by 
\[  \frac{1}{|\log|I||} \frac{H(P_{I}) +
K_{5}\de^{2}(P_{I},P_{J}) - H(P_{I}) -
K_{6}\de^{2}(P_{I},P_{J})}{2} \; .\]

It is evident that if we choose $K_{6}\leq K_{1}< K_{5}$,
estimate~\ref{equ:2,3} follows.       

\subparagraph{Proof of estimate ~\ref{equ:4,1}}
We need to show that
\[
\sum_{w\in P_{I}}|\log\mu(h(w))|\mu(w) - \sum_{w\in
P_{I}}|\log\mu(w)|\mu(w) \geq K_{5}\de^{2}(P_{I},P_{J})
\]

Here, we notice that it is  a well-known fact that the difference on
the left-hand side is non-negative. It can be checked directly by
calculus, or deduced from the variational principle for Gibbs measures
(see~\cite{bow}.)

Thus, we are trying to prove that this is a ``sharp'' inequality.

The idea is to split $P_{I}$ between two sets, called $E$ and $C$, so
that $h$ expands on $E$ and contracts on $C$. We define 
\[ E = \{w\in P_{I}\, : \, \frac{dh}{d\mu}(w) > 1\}\; \] then $C$ is
the complement of $E$. 

By Jensen's inequality
\[ \frac{\int_{E} \log\frac{dh}{d\mu}}{\mu(E)} \leq 
\log\frac{\mu(h(E))}{\mu(E)} \; .\] 

This allows an estimate of the average rate of expansion of $h$ on
$E$:

\begin{equation}\label{equ:5,1}
\frac{\mu(h(E))}{\mu(E)}\geq \exp\frac{\de(P_{I},P_{J})}{\mu(E)}
\; .
\end{equation}

Let us now look at the sum
\[ \sum_{w\in P_{I}}|\log\mu(h(w))|\mu(w) \; .\]

Its value given $P_{I}$ as well as sets $C,E,h(C),h(C)$ will be the
smallest if the Jacobian of $h$ is constant on both $A$ and $C$. Hence,
\[ \sum_{w\in P_{I}}|\log\mu(h(w))|\mu(w) - \sum_{w\in P_{I}}
|\log\mu(w)|\mu(w)\geq \]

\[ \mu(E)|\log\mu(h(E))| + 
(1-\mu(E))|\log(1-\mu(h(E)))| -\] 
\[- \mu(E)|\log\mu(E)| - (1-\mu(E))
|\log(1-\mu(E))|\; .\] 

To finish the proof of estimate ~\ref{equ:4,1}, we need to compare the
value of this difference (which must be non-negative) with 
$\de^{2}(P_{I},P_{J})$. 

Until the end of this proof we adopt notations
$x:=\mu(E)$ and $y:=\mu(h(E))$. We have $y>x$.
First of all, we see that
\[ x\log x + (1-x)\log (1-x) - x\log y - (1-x)\log (1-y) \geq
x(\frac{y}{x}-1-\log\frac{y}{x}) \]
provided that $y\geq x$. To see this, we notice that the equality
holds when $y=x$, and next we compare derivatives with respect to $y$.
As $x$ is fixed, the right-hand side of the preceding inequality grows
with $y/x$. This enables us to use estimate ~\ref{equ:5,1} and bound the
right-hand side of last inequality by 
\[ x\exp\frac{\de(P_{I}, P_{J})}{x} -x - \de(P_{I}, P_{J}) \; .\]    
As we neglect the terms of the exponential higher than the quadratic,
we get another estimate from below by
\[ \frac{\de^{2}(P_{I}, P_{J})}{2x} \] which is what was needed to prove 
estimate~\ref{equ:4,1}.

\subsection{ The upper exponent $\gu$.}

We begin with the observation that Fact~\ref{1} implies that the upper exponent
$\gu(x)$ is bounded from above by a uniform constant.
Here is the main result of this subsection. 
\begin{prop}
For almost all points of the circle the upper exponent
$\gu(x)$ is greater than $1$ and the estimate is uniform for maps
from our class.
\end{prop}

\paragraph{Checking procedure.}
Consider a sequence of nested partitions $\B(nr;f)$ and $\B(nr;\rho)$. 
Take an arbitrary interval $\Box_{f}^{nr}$ of the $nr$-th dynamical partition.
We will apply Proposition~\ref{prop:1,1} to partitions $\B((n+1)r;f)$
 and $\B((n+1)r;\rho)$ restricted to  $\Box^{nr}(f)$ and $\Box^{nr}(\rho)$
respectively. 

For rotations number of constant type Bounded Geometry implies that
the logarithms of conditional measures of atoms of our partitions are 
bounded  by a uniform constant. The same is with the logarithm of the
Jacobian of the isomorphism. 
So the hypothesis of Proposition~\ref{prop:1,1} is verified. 
We will keep the following scheme of {\em checking} the elements of
 the partitions $[\B((n+1)r;\rho):\Box^{nr}(\rho]$:
\begin{itemize} 
\item If the hypothesis of the implication in the thesis of 
Proposition~\ref{prop:1,1} is not satisfied for an element of
 $[\B((n+1)r;\rho):\Box^{nr}(\rho)]$  then we will call this element a ``good"
one. We stop {\em checking}.
\item Otherwise, we call an element of  $[\B((n+1)r;\rho):\Box^{nr}(\rho]$
a ``bad" one, denote by  $I^{(n+1)r}$, and pass to the subdivision of
 this interval by the next partition $\B((n+2)r;\rho)$. We repeat the
whole procedure.
\end{itemize}
Denote by $A$ a set of points which are covered infinitely many times by
``bad" elements of partitions $\B(nr;\rho)$.
\begin{lem}
The Lebesgue measure of $A$ must be zero.
\end{lem}
\begin{proof}
Suppose that the assertion of the Lemma is false. Then there is an arbitrary
fine cover of the set $A$ by ``bad" elements of the partition
 $\B(nr;\rho)$ (i.e. $n$ is large) which total length is greater then 
 $\lambda( A) > 0$.
We will apply Proposition~\ref{prop:1,1} step by step to the partitions 
$\B((n+1)r;\rho)$ restricted to elements $I^{nr}$. 
However, first we will make some preparation.

From Fact~\ref{1} it follows easily that 
\[\max_{I^{jr}\in \B(jr;\rho)} \frac{|\log|\phi^{-1}(I^{jr})||}
{|\log|I^{jr}||^{2}}\]
decreases up to a uniform constant as $1/j$. 

By the Discrepancy Lemma, 
\[ \de^{2}([\B((j+1)r;\rho): I^{jr})],
[\B((j+1)r;f):\phi^{-1}(I^{jr})]) \geq K. \]
Finally, repeated application of Proposition~\ref{prop:1,1} yields

\[\sum_{I^{nr}\in \B(nr;\rho)}
\frac{|\log|\phi^{-1}(I^{nr})||}{|\log|I^{nr}||}|I^{nr}|
 \geq \hspace*{\fill}\]
\[ K_{1} \lambda( A) \sum_{I^{1}\in B(1,r;\phi)} 
\frac{|\log|\phi^{-1}(I^{r})||}{|\log|I^{r}||} +  
 K_{2}\sum_{j=1}^{n-1}\frac{1}{j}.\]

The right-hand side of the above inequality tends to infinity
with $n$ while the left-hand side is bounded as we noticed at the beginning
of this subsection. This contradiction completes the proof.

\end{proof}

As a consequence, we see that the total length of ``good"
intervals of the partitions  $\B(nr;\rho)$ is equal to $1$.
Since now we will refer to ``good" intervals as  ``good" intervals of
 the first generation. We pass to a subdivision
of each ``good" interval of the first generation and repeat
the procedure of checking for all intervals of
the subdivision. By the same way as above we find ``good" intervals
of second generation which occupy again the whole space up to a set of
the Lebesgue measure zero. Repeating the procedure of {\em checking} countably
many times we will obtain a sequence of sets of ``good" intervals of different
generations. By the construction a ``good" interval of $n$-th generation
must be finer than any element of the partition $\B((n-1)r;\rho)$.

Denote by $G^{nr}_{x}$ a ``good" interval which belong to $\B(nr;\rho)$
and contains a point $x$ of the circle.
Let $B$ be a set of points which belong to infinitely many ``good"
intervals. Then for any $x \in$ and infinitely many $n$ 
Proposition~\ref{prop:1,1} implies the following estimate:
\[\gu(x) \geq \min(2, 1 + \frac{K}{H([\B((n+1)r;\rho): G^{nr}_{x}])}.\]
But the entropy $H([\B((n+1)r;\rho): G^{n}_{x}])$ is bounded
from above by a uniform constant. Hence,
 \[\gu(x) \geq 1+K_{2},\] 
where $K>0$ is an uniform constant.

\subsection{Lower exponent}
\paragraph{Statement.}
Now we are in a position to prove
\begin{prop}\label{prop:5,1}
For a constant $K>1$, we have
\[\gd(x) \geq K\]
for a full Lebesgue measure set of points $x$.
\end{prop}

\paragraph{Preliminaries of the proof.}
Since
\[ \de^{2}([\B((n+1)r;\rho: \Box^{nr}(\rho)],
[\B((n+1)r;f): \Box^{nr}(f)]) \geq K_{1} \, \]
and the entropy of $[\B ((n+1)r;\rho): \Box^{nr}(\rho)]$ is uniformly
bounded away from $0$, it follows that whenever 
\[ \frac{\log\phi^{-1}(|J|)}{\log|J|} < 1 + K_{2} \]
for uniform $K_{2}$, the assumptions of Proposition~\ref{prop:1,1} are
fulfilled for subpartions generated by $\B(n+1,r;\rho)$
and $\B(n+1,r;\rho)$ on $J$ and $\phi^{-1}(J)$ respectively.  
Now choose a number $a$ which is less than the a.e. upper exponent
and does not exceed $1+K_{2}$ either. Almost every trajectory will
spend an infinite amount of time above $a$.

Suppose that the lower exponent less or equal to $b-\epsilon$ on a positive
measure set $B$ ($\epsilon$ is arbitrary positive.) Our proof
will consist in showing that $b\geq a$. 

\subparagraph{The exponent as a random process.}
We define a random process
$(\tilde{Y}_{n})_{n=1,\ldots,\infty}$ so that each $\tilde{Y}_{n}$ is
measurable with
respect to $[\B(nr;\rho)$. If $J$ is an element of $\B((n+1)r;\rho)$,
$\tilde{Y}_{n}$ is constant on $J$ and equal to 
\[ \frac{\log\phi^{-1}(|J|)}{\log|J|} \; .\]
Then $\tilde{X}_{n}$ will be the increments of $\tilde{Y}_{n}$, i.e.
\[ \tilde{X}_{n} = \tilde{Y}_{n} - \tilde{Y}_{n-1}\; . \]  
We will use the following information about $\tilde{Y}_{n}$:
\begin{enumerate}
\item
$\tilde{X}_{n}$ is uniformly bounded by $K/n$. This follows immediately from
the definition of $\tilde{Y}_{n}$ and bounded geometry.
\item
$E(\tilde{X}_{n}|\tilde{Y}_{n-1}) \geq \frac{K'}{n}$ for a positive
 $K'$ provided that
$\tilde{Y}_{n-1}$ is less than $a$. This follows from
Proposition~\ref{prop:1,1}.
\end{enumerate}   

\subparagraph{The beginning of the proof.}
Suppose that $b<a$. 
Almost every trajectory of $(\tilde{Y}_{n})$ on $ B$ must
oscillate infinitely
many times between $a$ and $b$. Define an event $\tilde{A}_{n}$ as follows:
$\tilde{Y}_{n}>\frac{a+b}{2}$ and $\tilde{Y}_{n+1}\leq \frac{a+b}{2}$ and the
trajectory hits $b$ before hitting $a$. We will show that the series
of probabilities
\[ \sum_{n=1}^{\infty} P(\tilde{A}_{n})\] 
is summable which will immediately give us the desired contradiction. 

\paragraph{A supermartingale.}
We modify the process $(\tilde{Y}_{r})$ by making it constant after it hits
$a$ for the first time with $r\geq n$. The probability of the event
analogous to $\tilde{A}_{n}$ defined for the
new process will not change. To distinguish the modified process and
events defined for it from the old ones we will drop the tilde sign.  
Formally, $(Y_{r})_{r\geq n}$ also depends on $n$, but we choose not to
emphasize that in our notation.
The increments $X_{n}$ are still bounded by $K/n$, and $Y_{n}$ becomes
a submartingale (increasing conditional mean.)

\begin{defi}\label{defi:6,1}
We define a family of processes $(M(C,k,n))_{n=k,\ldots,\infty}$ indexed by $k$
by
\[ M(C,k,n) = \exp(\sqrt{n}(c-Y_{n}) -
kC\sum_{j=k+1}^{n}\frac{1}{j^{2}}) \; ,\]
where $c$ was used to denote $(a+b)/2$.
\end{defi}

\begin{lem}\label{lem:6,1}
One can choose uniform constants $K_{2}$ and $K_{3}$ so that for all
$k\geq K_{2}$ the process $M(K_{3},k,n)$ is a supermartingale.
\end{lem}
\begin{proof} 
We compute:
\[ E(M(C,k,n) | Y_{n-1}) =\]
\[ E(M(C,k,n-1))E(\exp(-\sqrt{k}X_{n})|Y_{n-1})\exp(-kC/n^{2}) \; .\]
One has to show that
\begin{equation}\label{equ:8,1}
 \log E(\exp(-\sqrt{k}X_{r})) \leq \frac{kC}{n^{2}} 
\end{equation}
if $k$ and $C$ are large. Since $X_{n}$ is of the order of $1/n \leq
1/k$, one can bound the exponent from above for large $k$ by
\[  1 - \sqrt{k}X_{n} + k(X_{n})^{2} \leq 1 - \sqrt{k}X_{n} +
kK_{3}/n^{2}\; .\]
Since $E(X_{n}|Y_{n-1})\geq 0$ we get
\[ E(\exp(-\sqrt{k}X_{n})) \leq 1 + kK_{3}/n^{2} \; .\]
Thus, whenever $k$ is large and $C \geq K_{3}$, Estimate~\ref{equ:8,1}
holds true, and the lemma immediately follows.  
\end{proof}

\paragraph{The bound for $P(A_{k})$.}
We substitute $A_{k}$ with a larger event $B_{k}$ which occurs when 
$\tilde{Y}_{k}\geq c$
and the trajectory by $(Y_{n})_{n\geq k}$ eventually hits $b$. 
We define the stopping time $j$ as the time of the
first crossing of $b$ by $Y_{n}$, $n>k$.
By the optional sampling theorem, (see~\cite{mar})
\[ \int_{B_{k}} M(C,k,j) \leq \int_{B_{k}} M(C,k,k)) \leq 1\]
since $M(C,k,k) \leq 1$ everywhere on $B_{k}$. One the other hand, 
\[M(C,k,j) \geq \exp(\sqrt{k}(c-b) -
nC\sum_{i=k+1}^{\infty}\frac{1}{i^{2}})\geq
\exp(\sqrt{k}(c-b) - C)\] 
on $B_{k}$. Thus, the measure of $B_{k}$ decreases like 
$K^{-\sqrt{k}}$ which is summable. 
\paragraph{}
In the consequence Proposition~\ref{prop:5,1} follows
and completes the proof of the {\bf Main Theorem}.
Changing the roles of $\rho$ and $f$ in the proof 
we immediately obtain the claim of the {\bf Remark}.

\section{ Hausdorff Dimension of $\mu$}
The H\H{o}lder type of singularity implies natural bounds on the Hausdorff
dimension of the measure $\mu$.
\begin{prop}\label{hausia}
The Hausdorff dimension of the invariant measure $\mu$ is equal to the 
lower exponent $\td(\mu)$ and, consequently, is uniformly bounded away
from $0$ and $1$.
\end{prop}
\begin{proof}
The proof the Proposition~\ref{hausia} is based on the following 
Frostman's Lemma:
\begin{fact}
Suppose that $\nu$ is a probabilistic Borel measure on the interval and for
$\nu$-a.e. $x$
\[\liminf_{\varepsilon \rightarrow 0} \log(\nu(x-\varepsilon,x+\varepsilon)/
\log(\varepsilon) = \kappa.\]  
Then the Hausdorff dimension of $\nu$ is equal to $\kappa$.
\end{fact}
By the Main Theorem follows that
\[\kappa = \td,\]
which completes the argument.
\end{proof}

\section{Open Questions.}

In the end of our presentation we would like to pose a few open questions
which we believe to be of natural interest and importance.

\begin{itemize}
\item Assuming that the rotation number is algebraic of degree $2$,
prove that the lower exponent is equal to its upper counterpart. This
should hold for the exponents related to $\lambda$ as well as $\mu$
and would give us just one exponent with respect to each measure.
\item In the same situation, establish a relation between exponents
 $\tau(\mu)$ and $\tau(\lambda)$.
\item Prove that $\tau(\mu)$ and $\tau(\lambda)$ are universal given
the rotation number (algebraic of degree $2$? any irrational?).
\item Do there exist critical circle homeomorphisms with a rotation
number of constant type
for which $\td(\mu) \neq \tg(\mu)$ and $\td(\lambda) \neq
\tg(\lambda)$? We suspect so. 
\item What is the situation for unbounded rotation numbers? Are main
results of this paper still valid? We suspect not. 
\end{itemize}

\end{document}